\documentclass{article}

\usepackage{indentfirst}
\usepackage{latexsym}
\usepackage{amsthm, amsmath, amsfonts, amssymb, graphicx}
\usepackage{algorithm, algorithmic}

\begin{document}

\title{Binary Fields on Limited Systems%
}


\author{ 
G\'abor P\'eter Nagy\\
\small{University of Szeged, Hungary}
\\
\small{\texttt{nagyg@math.u-szeged.hu}} \and
Valentino Lanzone \\
\small{University of Basilicata, Italy}
\\
\small{\texttt{valentinoweb@alice.it}} }

\date{}

\maketitle

\begin{abstract}
The intrinsic structure of binary fields poses a challenging complexity problem from both hardware and software point of view. Motivated by applications to modern cryptography, we describe some simple techniques aimed at performing computations over binary fields using systems with limited resources. This is particularly important when such computations must be carried out by means of very small and simple machines. The algorithms described in the present paper provide an increased efficiency in computations, when compared to the previously known algorithms for the arithmetic over prime fields.


\end{abstract}


\section{Introduction} \label{introd}
From the introduction of public key chryptography, numerous papers dealing with the problem of constructing efficient algorithms for the arithmetics of finite fields were published. With this respect, a vast amount of research has been carried out for Elliptic Curve Cryptography (ECC),~\cite{hkt}.

Recently, cryptosystems have been increasingly used in machines with very limited resources, like for instance smart cards, microchips and microcontrollers. This posed the problem of finding fast and efficient algorithms for field arithmetics when computations are to be performed by such simple devices.

The NIST\footnote{National Institute of Standards and Technology.} gave the recommendations for the selection of the underlying finite fields and elliptic curves. The latest revision of these standards was made available in the publication called FIPS 186-3~\cite{fips186-3}. This publication recommended 5 prime fields $\mathbb{F}_p$, with $p$ chosen among the following primes: $p_{192} = 2^{192}-2^{64}-1$, $p_{224} = 2^{224}-2^{96}+1$, $p_{256} = 2^{256}-2^{224}+2^{192}+2^{96}-1$, $p_{384} = 2^{384}-2^{128}-2^{96}+2^{32}-1$, $p_{521} = 2^{521}-1$, plus 5 binary fields: $\mathbb{F}_{2^{163}}$, $\mathbb{F}_{2^{233}}$, $\mathbb{F}_{2^{283}}$, $\mathbb{F}_{2^{409}}$ and $\mathbb{F}_{2^{571}}$. The NIST also gave detailed instructions on the use of elliptic curves over such finite fields.

Below we describe briefly some standard algorithms for the arithmetic of prime fields~\cite{bhlm}.

The primes $p$ for the prime fields are chosen with a bitsize divisible by 32. Further, $p$ must be either a Mersenne prime of the form $p = 2^{n}-1$, or a pseudo-Mersenne prime of the form $p = 2^{n}-r$ with the smallest possible integer $r$. We assume that the implementation platform has an $L$-bit architecture, with $L\in\{8,16,32,64\}$. Let $t = \lceil \log_{2}p \rceil$ and $m = \lceil t/L \rceil$, where $\lceil x \rceil$ denotes the least integer $k$ such that $k \geq x$; the elements of prime fields are the integers between $0$ and $p-1$ stored in software in an array of $m$ $L$-bit words: $a = (a_0, a_1, \ldots , a_{m-1})$. 

These primes allow an efficient modular reduction by using the replacement $a2^{n} \equiv ar \pmod{p}$, repeating it as necessary until the equivalent number modulo $p$ is obtained. 

Let $a = (a_0, a_1, \ldots , a_{m-1})$ and $b = (b_0, b_1, \ldots , b_{m-1})$ be two elements of a prime field $\mathbb{F}_p$. The addition is carried out by first finding the sum word by word and then reducing it modulo $p$. The modular addition is implemented by using the classic algorithm ``add with carry'', and the modular subtraction is implemented in a similar fashion where the carry is interpreted as a ``borrow''.

The multiplication is carried out by using the classic ``product term by term'', interpreted as ``product word by word'', and then reducing it modulo $p$.  We observe that, during the computation, we can easily represent each terms $a_i b_j=s_0+s_1 2^L$ still by the $L$-bit words $(s_0,s_1)$.

%

The inverse of a non zero field element $a \in \left\{  1, 2,\ldots ,p-1 \right\}$ is carried out by using a variant of the Extended Euclidean Algorithm. The algorithm maintains the invariants $Aa+dp=u$ and $Ca+ep=v$ for some $d$ and $e$ which are not explicitly computed. The algorithm terminates when $u=0$, in which case $v=1$, and $Ca+ep=1$, hence $C \equiv a^{-1}\pmod{p}$. Then, the division is carried out as $a/b=ab^{-1}$.

We have developed similar algorithms for binary fields in limited systems---whose small efficiency requires simple techniques---for the representation of bit sequences by suitable integers, with the property that addition and subtraction are the same, and with equality $1+1=0$.

In this paper we describe some simple algorithms that are designed to work with the arithmetic of the binary fields in limited systems such as microcontrollers, smart cards, etc. These algorithms are presented in form of pseudo-code.

\section{Arithmetic on binary fields and algorithms}
\label{Arith}

In a hardware circuit the data is represented by logical signals $\left\lbrace 0,1\right\rbrace $ and it uses the arithmetic of $t$-bits binary sequences. Therefore, the most appropriate choice for a finite field is $\mathrm{GF}(2^t)$. 
We have the following isomorphism:
\[\mathrm{GF}(2^t) \simeq \mathrm{GF}(2)[x]/p(x)\] where \[p(x) = x^t + r(x) = x^t + \sum_{i=0}^{t-1}p_ix^i, \quad p=(p_0,p_1,\dots,p_{t-1},1) \in \mathrm{GF}(2^{t+1}) \] is an irreducible polynomial of degree $t$ over $\mathrm{GF}(2)$. 
Using this isomorphism, the operations between $t$-bits binary sequences are identified with the operations between polynomials of degree $t-1$ modulo $p(x)$. 

To optimize the use of hardware memory, we can represent any sequence of $L$ bits with an unsigned integer between $0$ and $2^{L}-1$. More precisely, an element of $\mathrm{GF}(2^t)$ corresponds to $m=\lceil t/L \rceil$ unsigned integers. Then, using an appropriate representation of binary numbers as integers, we are able to access the bits representing the coefficients of the polynomials with appropriate functions and statements in terms of integers.

Let $d$ be the difference between $t$ and the degree of the polynomial $r(x)$. For practical reasons, polynomials $r(x)$ with few terms and degree as small as possible are preferable. One can use irreducible polynomials with three or five terms (trinomials and pentanomials, respectively) and such that $2d \geqslant t-1$.

The existence and the properties of certain irreducible polynomials, such as trinomials and pentanomials over $\mathrm{GF}(2)$, have been extensively investigated for at least 40 years following the paper of R.A. Shwan~\cite{Swan}. The relevant contributions prior to 1983 are surveyed in~\cite{ln}; see Chapter~3, Notes~5. Recent references on irreducible polynomials with few terms are~\cite{AM,BA,HPT,KK,PT}. In particular, a theorem due to Swan~\cite{Swan} implies that irreducible trinomials do not exist for  $t \equiv 0 \pmod{8}$. Furthermore, it follows from a result due to Bluher~\cite{BA} that they are rare when $t\equiv \pm\, 3 \pmod 8$; this fact originates from observations on trinomials and pentanomials arising from computations of Ahmadi and Menezes~\cite{AM}: If $t\equiv \pm\, 3 \pmod 8$  and $f(x)=\sum_{i=0}^t a_ix^i\in \mathrm{GF}[x]$ is an irreducible monic polynomial of degree $t$ such that $\mathrm{Tr}(a_i)=0$ for each $i$ with $1=i<t$, then $f$ contains a term $x^k$ with $t>k\geq t/3$ and $k=t-2 \pmod 4$. In particular, this shows for irreducible trinomials that the degree of the second term cannot be chosen to be of small.

When an irreducible trinomial of degree $t$ does not exist, the next best choice is a pentanomial. Usually, the polynomials are generated by deterministic irreducibility tests using computer computing, and a table of trinomials or pentanomials is available for $2 \leqslant t \leqslant 10000 $ in~\cite{hpl}.

We can write $r(x)=r_0x^{i_0}+r_1x^{i_1}+r_2x^{i_2}+r_3x^{i_3}$ with two zero terms in case of trinomials.

\subsection{Addition}
The addition of polynomials corresponds to the logical \textit{XOR} operation, also called \textit{exclusive or,} between bits of their corresponding binary sequences. Generally, programming languages for microcontrollers provide the \textit{XOR} operator for the integers. 

\begin{algorithm}
\caption{Addition in $\mathrm{GF}(2)[x]/p(x)$}
\label{alg1}
\begin{algorithmic}
\REQUIRE $ a=(a_0,a_1,\ldots,a_{m-1}),b=(b_0,b_1,\ldots,b_{m-1}),\quad a_i,b_i \in [0,2^{L}-1] $
\ENSURE $a+b=c=(c_0,c_1,\ldots,c_{m-1}),\quad c_i \in [0,2^{L}-1] $
\FOR{$i=0$ to $m-1$}
\STATE $c_i=a_i \mbox{\texttt{\^{ }}}  b_i$
\ENDFOR
\end{algorithmic}
\end{algorithm}

Algorithm \ref{alg1} computes the sum of two elements of $\mathrm{GF}(2^t)$ with computational complexity $O(m)$. The symbol ``\texttt{\^{ }}'' stands for the binary operator \textit{XOR} of unsigned integers.

\subsection{Reduction modulo $p(x)$}
Let $a(x) = \sum_{i=0}^{s}\alpha_ix^i$ be a polynomial of degree $s$, with $t \leqslant s \leqslant 2t-2$, represented by the binary sequence $(\alpha_0,\alpha_1,\ldots,\alpha_{s})$ with $\alpha_i \in \mathrm{GF}(2)$. 

Let $l=(\alpha_0,\alpha_1,\ldots,\alpha_{t-1})$ and $h=(\alpha_t,\alpha_{t+1},\ldots,\alpha_{2t-1})$, where $\alpha_i=0$ for $s+1 \leqslant i \leqslant 2t-1$, then we can write the polynomial as $a=l+hx^{t}$.

Since  $x^t \equiv r(x) \pmod{p(x)}$, we carry out the reduction of $a(x)$ modulo $p(x)$ using the following: 

\begin{itemize}
\item the equivalence $$ a \equiv l + hr(x) = l + r_0hx^{i_0} + r_1hx^{i_1} + r_2hx^{i_2} + r_3hx^{i_3} \pmod{p(x)} ;$$
\item the operations ``$f(x) \ll i$'' and ``$f(x) \gg i$'', which are the respective equivalents of shifting up and down $i$ positions in the binary sequence of the polynomial $f(x)$. 
\end{itemize}

\begin{algorithm}
\caption{Reduction modulo $p$ in $\mathrm{GF}(2)[x]$}
\label{alg2}
\begin{algorithmic}
\REQUIRE $ a=(\alpha_0,\alpha_1,\ldots,\alpha_{s}), \quad \alpha_i \in \mathrm{GF}(2), \quad s \leqslant 2t-2 $
\ENSURE $a \pmod{p(x)}$
\STATE $l=(\alpha_0,\alpha_1,\ldots,\alpha_{t-1})$ 
\STATE $h=(\alpha_t,\alpha_{t+1},\ldots,\alpha_{2t-1}), \mbox{ with } \alpha_i=0, \quad s+1 \leqslant i \leqslant 2t-1$
\WHILE {$ degree(a) \geqslant t $} 
\STATE $a=l$, $g=h \ll i_0$
\FOR {$i=0$ to $3$}
\IF {$r_i=1$} 
\STATE $a=a+g$
\ENDIF
\IF {$ i < 3$}
\STATE $g=g \ll (i_{i+1}-i_i)$
\ENDIF
\ENDFOR
\ENDWHILE
\end{algorithmic}
\end{algorithm}

When we shift a binary sequence by $i$ bits up or down, the ones into upmost or downmost $i$ bits, respectively, are lost. 
Our algorithms must guarantee that none of the ones are being shifted into oblivion, in order to assert that
\[[f(x) \cdot x^i]=[f(x) \ll i] \quad \mbox{and} \quad [f(x) / x^i]=[f(x) \gg i].\]

When a polynomial $a(x)$ has degree greater than $t-1$, we can delete the terms of degree greater than $t-1$ by using the equivalence $ a \equiv l + hr(x) \pmod{p(x)} $ and repeating it if necessary. Since $2d \geqslant t-1$, we need to iterate this operation no more than twice. So, we obtain Algorithm~\ref{alg2}, which has computational complexity $O(km)$, with $k=4$ or $k=8$ according as $p(x)$ is a trinomial or a pentanomial. 

\subsection{Square}
Since $\mathrm{GF}(2)$ is a field of characteristic 2, the following equality holds  
\[ \left(\sum_{i=0}^{t-1}{\alpha_{i}x^{i}}\right)^2 = \sum_{i=0}^{t-1}{\alpha_{i}x^{2i}}, \quad \alpha_{i} \in \mathrm{GF}(2). \] 

\begin{algorithm}
\caption{Square in $\mathrm{GF}(2)[x]/p(x)$}
\label{alg3}
\begin{algorithmic}
\REQUIRE $ a=(\alpha_0,\alpha_1,\ldots,\alpha_{t-1}), \quad \alpha_i \in
\mathrm{GF}(2) $
\ENSURE $a^2 \pmod{p(x)}$
\STATE temporary variable: $b=(\beta_0,\beta_1,\ldots,\beta_{2t-2}), \quad
\beta_i \in \mathrm{GF}(2) $
\FOR {$i=0$ to $t-2$}
\STATE $ \beta_{2i} = \alpha_i $
\STATE $ \beta_{2i+1} = 0 $
\ENDFOR
\STATE $ \beta_{2t-2} = \alpha_{t-1} $
\STATE $a^2 = b \pmod{p(x)}$
\end{algorithmic}
\end{algorithm}

Therefore, we can compute the square of a polynomial simply by doubling its indices and then performing the reduction modulo $p(x)$. 
We obtain Algorithm~\ref{alg3}, whose main computational cost is due to reduction. 

\subsection{Product}
Let $ a,b \in \mathrm{GF}(2)[x] $ be two polynomials, with 
\[a=\sum_{i=0}^{t-1}{\alpha_{i}x^{i}}. \] 

Since $[b \cdot x^i]=[b \ll i]$ the product between $a$ and $b$ is 
\[a \cdot b = \sum_{i=0}^{t-1}{\alpha_{i}x^{i}} \cdot b =
\sum_{i=0}^{t-1}{\alpha_{i}} \cdot (bx^{i}) = \sum_{i=0}^{t-1}{\alpha_{i}} \cdot
(b \ll i),\] 
which has computational complexity $O(t)$ plus $t$ shifts and the reduction's cost.

But, we can perform the product faster as follows. Let $ a=(\alpha_0,\alpha_1,\ldots,\alpha_{t-1}) \in \mathrm{GF}(2)^t$, $b=(b_0,b_1,\ldots,b_{m-1}) \in [0,2^{L}-1]^m $, $w=\lceil 2t/L \rceil$, and define the operation $ Sh(b,i)= b \ll iL $.
We note that $ Sh(b,i)=(s_0,s_1,\ldots,s_{w-1})$, where $s_j=b_{j-i}$, if $j \in [i,m-1+i]$ and $s_j=0$ otherwise. 

By using the operation $Sh$, we only need to do $L$ shift operations, instead of $t$, in this way:
$$
a \cdot b = \sum_{e=0}^{L-1}{ \left\lbrace  \sum_{i=0}^{m-1}{\alpha_{iL+e} \cdot \left[  Sh\big((b \ll e),i\big) \right]  } \right\rbrace  }
$$

\begin{algorithm}
\caption{Product in $\mathrm{GF}(2)[x]/p(x)$}
\label{alg4}
\begin{algorithmic}
\REQUIRE $ a=(\alpha_0,\alpha_1,\ldots,\alpha_{t-1}) \in \mathrm{GF}(2^t), \ b=(b_0,b_1,\ldots,b_{m-1}) \in [0,2^{L}-1]^m $
\ENSURE $a \cdot b \pmod{p(x)}$
\STATE temporary variables: $c, d \in [0,2^{L}-1]^w$
\STATE $ c=(c_0,c_1,\ldots,c_{w-1}), \mbox{ with } c_i=0, \ 0 \leqslant i \leqslant w-1 $
\STATE $ d=(d_0,d_1,\ldots,d_{w-1}), \mbox{ with } d_i = b_i, \ 0 \leqslant i \leqslant m-1, \mbox{ and } d_i = 0, \ m \leqslant i \leqslant w-1$
\FOR {$e=0$ to $L-1$}
	\FOR {$i=0$ to $m-1$}
		\IF {$\alpha_{iL+e}=1$}
			\FOR {$j=i$ to $w-1$}
			\STATE $c_j=c_j+d_{j-i}$
			\ENDFOR
		\ENDIF
	\ENDFOR
	\STATE $d=d \ll 1$
\ENDFOR
\STATE $a \cdot b = c \pmod{p(x)}$
\end{algorithmic}
\end{algorithm}

We have Algorithm~\ref{alg4}, which has computational complexity $O(Lm) \sim O(t)$ plus $L$ shifts and the reduction's cost.


\subsection{Inversion and division}
To compute the inverse of polynomials we use a variant of the classical
Euclidean algorithm. We can carry out the division between two polynomials by multiplying the first one by the inverse of the second one.

Let $a(x)$ and $b(x)$ be two polynomials in
$\mathrm{GF}(2)[x]$. Then, $\gcd (a,b)=\gcd (b-ca,a)$ for all polynomials $c \in
\mathrm{GF}(2)[x]$. If $\deg (b)\geqslant \deg (a)$ and $j=\deg (b)-\deg (a)$,
we can compute $r=b+x^ja$ and hold $\gcd (a,b)=\gcd (r,a)$. 

\begin{algorithm}
\caption{Inversion in $\mathrm{GF}(2)[x]/p(x)$}
\label{alg5}
\begin{algorithmic}
\REQUIRE $ a=(\alpha_0,\alpha_1,\ldots,\alpha_{t-1}) \neq 0, \quad \alpha_i \in
\mathrm{GF}(2)$
\ENSURE $a^{-1} \pmod{p(x)}$
\STATE temporary variables in $\mathrm{GF}(2^t)$: $u=a, \quad v=p, \quad g_1=1,
\quad g_2=0 $
\WHILE {$ degree(u) \neq 0 $}
\STATE $j=degree(u)-degree(v)$
\IF {$j<0$}
\STATE $swap(u,v)$, $swap(g_1,g_2)$, $j=-j$
\ENDIF
\STATE $u=u+(v \ll j)$
\STATE $g_1=g_1+(g_2 \ll j)$
\ENDWHILE
\STATE $a^{-1}=g_1$
\end{algorithmic}
\end{algorithm}

With this variant, we can use the extended Euclidean algorithm and obtain Algorithm~\ref{alg5}
which has computational complexity $O(4tm)$, see~\cite{GuideECC}. 

\section{Tests performed}
We tested these algorithms on a commercially available and very cheap board. 
Such a board, called Arduino\texttrademark \ Duemilanove\footnote{http://www.arduino.cc/}, 
has computing power similar to smart cards and has the following features:  
\begin{itemize} 
  \item ATmega168 microcontroller\footnote{Low Power AVR\textregistered \
    Microcontroller manufactured by ATMEL\textregistered.} ; 
  \item 16 KB (available 14 KB) in system self-programmable flash memory;
  \item 1 KB SRAM and 512 Bytes EEPROM;
  \item 16 MHz clock speed;
  \item language based on C/C++;
  \item standard serial communication.
 \end{itemize}


Below, we show the most significant results obtained on the 5 binary fields that NIST recommended in the publication FIPS 186-3, with following polynomial basis representation:

\begin{itemize}
\item $\mathbb{F}_{2^{163}} = \mathrm{GF}(2)[x]/(x^{163}+x^{7}+x^{6}+x^{3}+1)$,
\item $\mathbb{F}_{2^{233}} = \mathrm{GF}(2)[x]/(x^{233}+x^{74}+1)$,
\item $\mathbb{F}_{2^{283}} = \mathrm{GF}(2)[x]/(x^{283}+x^{12}+x^7+x^5+1)$,
\item $\mathbb{F}_{2^{409}} = \mathrm{GF}(2)[x]/(x^{409}+x^{87}+1)$,
\item $\mathbb{F}_{2^{571}} = \mathrm{GF}(2)[x]/(x^{571}+x^{10}+x^5+x^2+1)$.
\end{itemize}

\begin{table}[ht]
\begin{center}\begin{tabular}{|l|c|c|c|c|c|}
\hline
\textit{Degree of field}     & \textit{163} & \textit{233} & \textit{283} & \textit{409} & \textit{571} \\
\hline
multiplication on binary fields       & 16  & 29  & 40  & 80  & 149 \\
inversion on binary fields         & 60  & 105 & 145 & 282 & 505 \\
\hline
\end{tabular}\end{center}
\caption{Execution times on binary fields (in ms)} 
\label{time-binaryfields}
\end{table}

\begin{table}[ht]
\begin{center}\begin{tabular}{|l|c|c|c|c|c|}
\hline
\textit{Degree of field}     & \textit{192} & \textit{224} & \textit{256} & \textit{384} & \textit{521}\\
\hline
multiplication on prime fields        & 6   & 7   & 9   & 18   & 29   \\
inversion on prime fields          & 234 & 344 & 490 & 1442 & 3258 \\
\hline
\end{tabular}\end{center}
\caption{Execution times on prime fields (in ms)} 
\label{time-primefields}
\end{table}

\begin{figure}[htbp]
\begin{center}
\scalebox{0.70}{\includegraphics[bb=25 530 440 795]{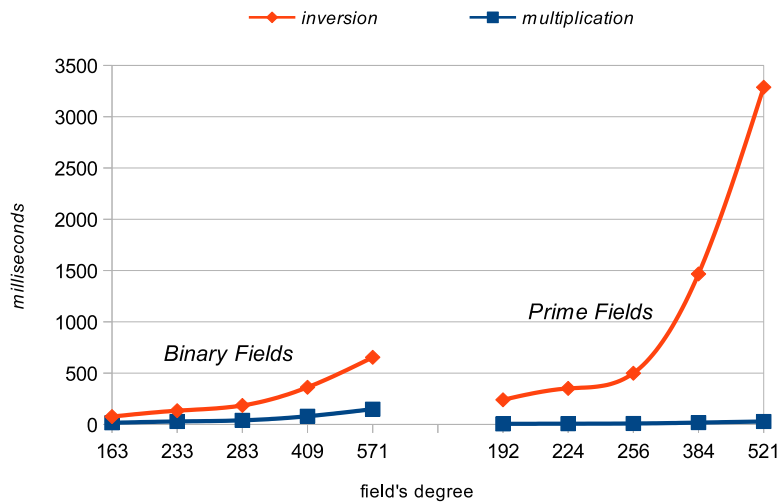}}
\caption{Time comparisons, $L=16$ bits.}
\label{fig:graph}
\end{center}
\end{figure}

In order to do a comparison, we have also implemented the algorithms on the NIST prime fields shown in the Introduction~\ref{introd}.
In Tables~\ref{time-binaryfields} and~\ref{time-primefields}, we put the execution times to multiply and invert on the NIST binary fields and on the NIST prime fields respectively.
In Figure~\ref{fig:graph} we provide a visual comparison between the execution times on binary fields and prime fields.

\section{Conclusion}
In this paper, we presented an implementation of the arithmetic in $\mathrm{GF}(2^t)$ with basic polynomial, using straightforward algorithms with low use of memory. 

The algorithms we used are as generic as possible, so we can easily change the parameters and the underlying field $\mathrm{GF}(2^t)$.  For their flexibility, these algorithms can be used in systems with limited computing resources. 

From the comparison between the execution times, we observe that the multiplication on prime fields requires an execution time which is shorter than on binary fields, while the operation of inversion on prime fields has an execution time much larger than on binary fields, and this grows very rapidly.

Furthermore, we can observe that our algorithms proved to be very efficient and particularly suitable for small devices and tasks which require the use of arithmetic inversions.

\section*{Acknowledgements}

The authors would like to thank the referee for many helpful comments and hints. This work was finantially supported by the T\'AMOP-4.2.2/08/1/2008-0008 program of the Hungarian National Development Agency, by the Italian Ministry MIUR PRIN 2011-12, by Italian Institute of high Mathematics INdAM-GNSAGA. The paper is based on the talk given by \mbox{V.~Lanzone} in the The Second Conference of PhD Students in Mathematics (CSM 2), Szeged, Hungary, 2012.


\end{document}